\documentclass[12pt,letterpaper]{amsart}
\usepackage{amssymb,latexsym, amsmath, amsxtra}
\usepackage[dvips]{graphics}

\theoremstyle{plain}

\theoremstyle{definition}

\theoremstyle{remark}

\numberwithin{equation}{section}
\numberwithin{figure}{section}

\newcommand{\Q}{\mathbb Q}

\newcommand{\N}{\mathbb N}

\newcommand{\Fa}{\mathcal F_1}
\newcommand{\Fb}{\mathcal F_2}
\newcommand{\Fc}{\mathcal F_3}
\newcommand{\Fd}{\mathcal F_4}

\def\({\left(}
\def\){\right)}

\def\yr #1#2#3#4{(#1#2#3#4)}
\def\publ{\relax}
\def\publaddr{\relax}
\def\vol{\relax}

\def\jour{\relax}

\newcommand{\ontop}[2]{\genfrac{}{}{0pt}{}{#1}{#2}}

\begin{document}
\title{Modeling families of L-functions}
\author{David W Farmer}

\thanks{Work supported by the
American Institute of Mathematics
and by the Focused Research Group grant (0244660) from the NSF.  
This paper is an expanded version of a talk given at the workshop
``Random Matrix Theory and Elliptic Curves''
held at the
Isaac Newton Institute, February 2004.
}

\thispagestyle{empty}
\vspace{.5cm}
\begin{abstract}
We discuss the idea of a ``family of $L$-functions'' and
describe various methods which have been used to make
predictions about $L$-function families.  The methods
involve a mixture of random matrix theory and heuristics
from number theory.  Particular attention is paid to 
families of elliptic curve $L$-functions.
We describe two random matrix models for
elliptic curve families:
the Independent Model  and the Interaction Model.
\end{abstract}

\address{
{\parskip 0pt
American Institute of Mathematics\endgraf
360 Portage Ave.\endgraf
Palo Alto, CA 94306\endgraf
farmer@aimath.org\endgraf
}
  }

\maketitle

\section{Introduction}\label{sec:intro}

Using ensembles of random matrices to model
the statistical properties of a family of $L$-functions
has led to a wealth of interesting conjectures and results
in number theory.
In this paper we survey recent results in the hopes of 
conveying our best current answers to these questions:
\begin{enumerate}
\item What is a family of $L$-functions?
\item How do we model a family of $L$-functions?
\item What properties of the family can the model predict? \label{q:properties}
\end{enumerate}

In the remainder of this section we briefly review some commonly studied
families and describe some of the properties which have 
been modeled using ideas from random matrix theory.  
In Section~\ref{sec:family} we provide a definition of 
``family of $L$-functions'' which has been 
successful in permitting precise conjectures, and we briefly
describe how to model such a family.
In Section~\ref{sec:EC} we discuss families of elliptic curve
$L$-functions and show that there is an additional subtlety
which requires us to slightly broaden the class of random matrix
models we use.  Then in Section~\ref{sec:refined} we discuss
how to go beyond the leading-order terms which random matrix
theory can model, and how one can avoid using
random matrix theory when modeling a family of $L$-functions.

If you only care about elliptic curves and their $L$-functions 
you can safely skip to Section~\ref{sec:EC}.

I thank Brian Conrey, Nina Snaith, Matt Young, and
Steven J.~Miller for many helpful conversations.

\subsection{A quick history of families}

The idea that collectively the zeros of a single $L$-function behave in 
a manner that can be modeled statistically started with 
Montgomery's~\cite{Mon} work on the pair correlation of
zeros of the zeta-function.  Combined with the large-scale
numerical calculations
of Odlyzko~\cite{Odl}, this provided convincing evidence that, 
to leading order, the local statistics of the zeros
of the zeta-function, suitably rescaled, were the same as
those of large random unitary matrices. 

A similar collective behavior was noted long ago 
and is termed the ``$q$-analogue'' for Dirichlet
$L$-functions.  That is, results about all Dirichlet
$L$-functions $\mod q$ look just like results for
the Riemann $\zeta$-function in $t$-aspect. 
For example, the formulas for moments of $|\zeta(\frac12+it)|$
and $|L(s,\chi)|$ are identical 
but for replacing $t$ by~$q$.  
Another early example of collective behavior is the
pair correlation of zeros of quadratic Dirichlet $L$-functions~\cite{OS1}.
Clearly something interesting is going on. 

The idea of a family of $L$-functions with an associated symmetry type
began with the work of Katz and Sarnak~\cite{KSa}.  They consider 
families of function field $L$-functions, where in this case a
``family'' of $L$-functions is the set of $L$-functions associated to 
a set of curves having certain properties.
Here the collection of curves must be ``natural'' in the sense
that the {monodromy group} of
the family ties it all together.
They show that 
to leading order
the statistics of the (normalized) zeros of these $L$-functions,
when averaged over the family, are the same as the statistics of the
(normalized) eigenvalues of random matrices chosen from a 
classical compact group.
Here the matrices are chosen uniformly with respect to Haar measure,
and the size of the matrices scales with the conductor of the
$L$-function.

For global $L$-functions there does not (yet?) exist an analogue
of ``monodromy'', but it still has been found that to naturally
occurring families one can associate a classical group of matrices. 
The zeros of the $L$-functions have, to leading order, the same
statistics as the eigenvalues of a randomly chosen matrix from the group.
And the appropriately rescaled critical values of the $L$-functions
have, to leading order, the same distribution as the ``critical values''
of the characteristic polynomials of the matrices from the group,
chosen uniformly with respect to Haar measure.  The correspondence between
the matrices and the $L$-functions involves equating the eigenvalue spacing 
with the zero spacing, or equivalently setting the matrix size equal
to the conductor of the $L$-function.  See Section~\ref{ssec:modeling}
for more details.

After the work of Katz and Sarnak there quickly appeared many examples
of $L$-functions families behaving in a manner predicted by
random matrix theory.  Some of the families considered were:
$L$-functions associated to holomorphic cusp forms 
(in either weight or level aspect);
Dirichlet $L$-functions (either all or quadratic);
and various twists or symmetric powers of $L$-functions.
Low-lying zeros were considered by 
Iwaniec, Luo, and Sarnak~\cite{ILS}, 
Rubinstein~\cite{Rub},
\"Ozl\"uk and Snyder~\cite{OS2}, and others.
Moments were considered by Iwaniec and Sarnak~\cite{IS},
Kowalski, Michel, and VanderKam~\cite{KMV}, 
Soundararajan~\cite{S}, and others.  

There is an important distinction between the predictions for
zeros of $L$-functions and the predictions for moments.
For the zeros there is a natural way to normalize:
rescale so that the average spacing is~$1$.  This rescaling involves
the conductor of the $L$-function and the degree of
the characteristic polynomial, and this is the source of the principle
that one chooses the size of the random matrix to equal the
(logarithm of the) conductor of the $L$-function. 
With these normalizations one obtains accurate predictions for the
leading-order behavior of statistics of zeros of $L$-functions.
However, in addition to its zeros a polynomial is determined
by an overall scale factor. So one may hope to use the characteristic
polynomials to model the $L$-functions, but there will be a
correction factor that does not come from random matrix theory.
The use of the characteristic polynomial to model the
$L$-function was begun by Keating and Snaith~\cite{KS1}.
The first situation which was well understood 
is the
moments of $L$-functions, for which there are
explicit predictions for the arithmetic scale factor.
See~\cite{CG, CF, KS2}.

\section{What is a family?}\label{sec:family}

There does not yet exist an adequate definition of ``family of $L$-functions''.
An attempt is made in~\cite{CFKRS} to define a family axiomatically,
and we will describe that definition here.
In that definition the axioms are chosen so that it is possible to produce a plausible
conjecture for the critical moments of the family.  
%We will describe the
%case of moments of elliptic curve $L$-functions below.

Complete details and many examples are in~\cite{CFKRS}, so we will
just highlight the key features of that definition of a family.
The main idea is that one starts with a fixed $L$-function and a 
family of ``characters'', and the family of $L$-functions is produced
by twisting the fixed $L$-function by the family of characters.
Note that here the term \emph{character} is used to cover more
general classes of functions than just Dirichlet characters.

\subsection{$L$-functions} 
We wish to define a ``family of $L$-functions'', so first we have to
give the definition we will use for ``$L$-function''.  
The definition of an $L$-function which we give below is 
slightly different than what is known as
the ``Selberg class,'' but it is conjectured that the
two are in fact equal.

Let $s=\sigma+it$ with
$\sigma$ and $t$ real. An \emph{$L$-function} is a Dirichlet
series
\begin{equation}
L(s)=\sum_{n=1}^\infty \frac{a_n}{n^s},
\end{equation}
with $a_n =O_\varepsilon (n^\varepsilon)$ for every $\varepsilon>0$, which
has three additional properties.

\noindent\emph{Analytic continuation:}  $L(s)$ continues to a meromorphic
function of finite order with at most finitely many poles, and all
poles are located on the $\sigma=1$ line.

\noindent\emph{Functional equation:} There is a number $\varepsilon$
with $|\varepsilon|=1$, and a function $\gamma_L(s)$ of the form
\begin{equation}
\gamma_L(s)= P(s) Q^s \prod_{j=1}^w \Gamma(\tfrac12 s+ \mu_j) ,
\label{eqn:gammafactors}
\end{equation}
where $Q>0$, $\Re \mu_j \ge 0$, and $P$ is a
polynomial whose only zeros in $\sigma>0$ are at the poles of
$L(s)$, such that
\begin{equation}
\xi_L(s):= \gamma_L(s)L(s)
\end{equation}
is entire, and
\begin{equation}
\xi_L(s)    =\varepsilon \overline{\xi_L}(1-s),
\end{equation}
where $\overline{\xi_L}(s) := \overline{\xi_L (\overline{s})}$ and
$\overline{s}$ denotes the complex conjugate of~$s$ .

The number $w$ is called the \emph{degree} of the $L$-function.
That number will also appear in the Euler product.

\noindent\emph{Euler product:}  For $\sigma>1$ we have
\begin{equation}
L(s)=\prod_p \prod_{j=1}^w (1-\gamma_{p,j}p^{-s})^{-1} ,
\end{equation}
where the product is over the primes~$p$, and
each
$|\gamma_{p,j}|$ equals~$1$ or~$0$.

Note that $L(s)\equiv 1$ is the only constant $L$-function,
 the set of $L$-functions is closed under products,
and if $L(s)$ is an $L$-function then so is $L(s+iy)$ for
any real~$y$.  An $L$-function is called \emph{primitive}
if it cannot be written as a nontrivial product of
$L$-functions.  Throughout this paper we assume all $L$-functions
are primitive, although we usually omit the word ``primitive.''

\noindent\emph{Conductor:} Associated to an $L$-function is
its \emph{conductor}, a number which measures the ``size'' of the
$L$-function.  The paper~\cite{CFKRS} introduced a refined notion
of conductor which, to leading order, is the logarithm of the
usual notion of conductor.  The refined conductor is necessary in 
order to have any hope of conjecturing the full main term in a
general mean
value of the $L$-function.
Write the functional equation in
asymmetric form:
\begin{equation}
L(s)=\varepsilon X_L(s) \overline{L}(1-s),
\end{equation}
where $ \displaystyle X_L(s) =
\frac{\overline{\gamma_L}(1-s)}{\gamma_L(s)}$.
Then the refined conductor of $L(s)$, denoted $c(L)$, is given by
$c(L)=|X_L'(\tfrac12)|$.

\subsection{Families of characters}
By a family of characters we mean a collection
of arithmetic functions $\mathcal F$, where each $f\in \mathcal F$ is
a sequence $f(1)=1$, $f(2)=a_{2,f}$, $f(3)=a_{3,f},\ldots$ 
whose generating function
\begin{equation}
L_f(s)=\sum_{n=1}^\infty \frac{a_{n,f}}{n^s}
=
\prod_p \prod_{j=1}^v (1-\beta_{p,j} p^{-s})^{-1} 
\end{equation}
is a (primitive) $L$-function such that the collection
$\{L_f \ : \ f\in \mathcal F\}$ has some nice properties.
If we order the $L$-functions $L_f$ by conductor~$c(f)$,
then the data $\{Q;\, \mu_1,\ldots,\mu_w\}$ 
in the functional equation of $L_f$ should be 
monotonic functions of the conductor, and  
the counting function $M(X):=\#\{f\in \mathcal F\ |\ c(f)\le X\}$ 
should be nice.  The final condition on the family of characters
is the existence of an orthogonality relation among the $f\in \mathcal F$.
Specifically, we require that
if $m_1,\dots,m_k$ are integers then
the average
\begin{equation}\label{eqn:orthogonality}
\delta_{\ell}(m_1,\dots,m_k)
:=
\lim_{X\to \infty}
M(X)^{-1}\sum_{\ontop{f\in \mathcal{F}}{c(f)\le X}}
f(m_1)\dots
f(m_\ell)\overline{f(m_{\ell+1})\dots f(m_k)} 
\end{equation}
exist and be multiplicative.  That is,
if $(m_1 m_2 \dots m_k, n_1 n_2 \dots n_k)=1$, then
\begin{equation}\label{eqn:deltamultiplicative}
\delta_\ell(m_1n_1,m_2n_2,\dots, m_kn_k)=
\delta_\ell(m_1,\dots,m_k)\delta_\ell(n_1,\dots, n_k).
\end{equation}
See Section~3.1 of~\cite{CFKRS} for more details.
%Note that some families are unions of increasingly large 
%pieces having the same conductor.  Examples are the Dirichlet 
%$L$-functions and the $L$-functions associated to holomorphic
%cusp forms (in either weight or level aspect).  
%For those families it is belived that we can
%modify the orthogonality relation~\eqref{eqn:orthogonality}
%so that $X$ runs over the possible conductors, and the 
%sum is restricted to the $L$-functions of conductor~$X$.

\subsection{Families of $L$-functions}
Now we create a family of $L$-functions by starting with
a fixed $L$-function
\begin{equation}
L_g(s)=\sum_{n=1}^\infty \frac{a_{n,g}}{n^s} =
\prod_p \prod_{j=1}^w (1-\gamma_{p,j}p^{-s})^{-1}.
\end{equation}
Then the elements of our $L$-function family
$\mathcal L(\mathcal F)$ are the Rankin-Selberg convolutions
\begin{align} \label{eqn:rankinselberg}
\mathcal{L}(s,f)=
L_{f\times g}(s)
=&
\prod_p\prod_{i=1}^v\prod_{j=1}^w(1-\beta_{p,i}\gamma_{p,j} p^{-s})^{-1}\cr
=&\sum_{n=1}^\infty
\frac{a_{n,f\times g}}{n^s}.
\end{align}
(There may be some issues with the local factors at the bad primes).
Note that if $w=1$ or $v=1$ then 
\begin{equation}
\mathcal{L}(s,f)=\sum_{n=1}^\infty
\frac{a_{n,f} a_{n,g}}{n^s}.
\end{equation}
And in particular if $L_g$ is the Riemann zeta-function, then
$\mathcal{L}(s,f)=L_f(s)$.

The point of this definition of ``family'' is that the axioms 
provide the necessary ingredients to apply the recipe 
in~\cite{CFKRS} to conjecture the full main term in the
shifted $K$th moment
\begin{equation}
{M(X)^{-1}} \sum_{c(f)\le X} \prod_{1\le k\le K}
L(\tfrac12+\alpha_k) ,
\end{equation}
or more generally a 
shifted ratio~\cite{CFZ}
\begin{equation}
{M(X)^{-1}} \sum_{c(f)\le X} \prod_{1\le k\le K}
\frac{L(\tfrac12+\alpha_k)}{L(\tfrac12+\delta_k)} .
\end{equation}
Having such a mean value is sufficient to conjecture just about
anything you would like to know about the zeros and the value
distribution of the $L$-function.  See~\cite{CS} for examples.

Note that some families are unions of increasingly large
pieces having the same conductor.  Examples are the Dirichlet
$L$-functions and the $L$-functions associated to holomorphic
cusp forms (in either weight or level aspect). For those families 
it is believed that the heuristics for moments will produce
a reasonable conjecture for the average over a fixed (large)
conductor.

Although this definition of ``family'' is useful for certain
applications, it lacks the concreteness of the function field case.
In particular, there does not yet exist an analogue of monodromy for
such a family, and computing the symmetry type of the family
is not straightforward.  We discuss this in the next section.

\subsection{Modeling a family of $L$-functions}\label{ssec:modeling}
Given a family of $L$-functions one can ask questions about its value
distribution or about the distribution of its zeros.  In most
cases current technology is not sufficient to answer 
the interesting questions, so the next hope is to find a
plausible conjecture.  Only recently have such conjectures been
found, and the new ingredient is to use random matrices to model
the family of $L$-functions.

The idea is to associate a classical compact group,
$U(N)$, $Sp(2N)$, $O(N)$, $SO(2N)$, or $SO(2N+1)$,
 to the family.  The local statistics of the eigenvalues
should agree, to leading order, with the corresponding
local statistics of the zeros of the $L$-functions.
And, to leading order and after compensating by an arithmetic
constant, the value distribution of the characteristic
polynomial
\begin{equation}\Lambda (z)=\Lambda_A(z)=\det(I-A^*z)
=\prod_{n=1}^N \(1-z e^{-i\theta_n}\)
\end{equation}
near the point $z=1$ should agree with the value distribution
of the $L$-functions near the critical point.
Here $A$ is an $N\times N$ unitary matrix $A$ 
and  $A^{*}$ is the Hermitian
conjugate of~$A$,
so the eigenvalues of $A$ lie on the unit circle
and are denoted by~$e^{i\theta_n}$.

In the above correspondence the size of the matrix is set
equal to the conductor of the $L$-function.  (Actually, to an integer close
to the conductor, but to leading order such discrepancies do not matter).
To see why this is a natural choice, 
consider the functional equation satisfied by the characteristic polynomial:
\begin{equation}
\Lambda_A(z) = (-1)^N \det(A) z^N \Lambda_{A^*}(z^{-1}). 
\end{equation}
If we identify $(-1)^N \det(A)$ with $\varepsilon$ and $z^N$ with
$X_L(s)$ then we have a perfect correspondence between the 
functional equations of $\Lambda_A(z)$ and $L(s)$, the unit circle
playing the role of the critical line and $z=1$ the critical point.
Just as for $L$-functions, we define the conductor as 
$\frac{d}{dz} z^N$ evaluated at the critical point, so $N$
is the conductor.  Note that identifying conductors is equivalent
to equating the average spacing between the zeros.
Values near the critical point are modeled using the 
correspondence $\Lambda(e^{-z}) \leftrightarrow L(\frac12+z)$.

It remains to identify the matrix group which corresponds to the family.
From the functional equation it is almost possible to determine
the group: the only ambiguity is to distinguish between
$SO(2N)$ and $Sp(2N)$.  At one time it was thought that this
case could be easily resolved because $SO(2N)$ families always arise
as ``half'' of a larger family, the other half being modeled 
by $SO(2N+1)$. On the other hand $Sp(2N)$ families do
not have such a ``partner''.  A counterexample to that hope is
described in~\cite{MD2}.  But even if that approach were viable,
it is unsatisfactory because it relies on the fact that
the symmetry type can be found among a small list of possibilities.
Fortunately, there are other methods.
One possibility is to compute the 1-~and 2-level densities of the family.
This usually can be done rigorously for functions with small support,
and this is sufficient to distinguish among the classical compact groups.
But again we are relying on the fact that the symmetry type can
be found on a short list.  That objection can be overcome if one can 
conjecture the level densities in the full range, but that can
be quite difficult in practice.
Another possibility is to use the
recipe in~\cite{CFKRS} to conjecture the moments of the family.
This unambiguously identifies the group, and it also can tell you if
the family is not modeled by one of those groups.  
Unfortunately, it is not clear that the recipe in~\cite{CFKRS}
can be applied to all interesting families of elliptic curve $L$-functions,
such as the family~$\Fc$ given in~\eqref{eqn:allrankr}.

\subsection{Summary of modeling}

Just to be pedantic, we note the following answers to the 
questions posed at the beginning of Section~\ref{sec:intro}:

\begin{enumerate}
\item A family of $L$-functions is a set of $L$-functions,
ordered by conductor, which is built
in a particular way from a family of characters. The counting
function of the family should be nice, and the data in the functional
equation should be monotonic functions of the conductor.

\item A family is modeled by associating to it a classical compact
matrix group.  The specific compact group can usually be determined by computing
the level densities of the low-lying zeros of the family, or by
conjecturing the moments of the family.   The size of the matrices scales with the
(logarithmic) conductor of the $L$-functions.

\item  To leading order the rescaled local zero statistics  of the family
are the same as the rescaled local eigenvalue statistics of the 
group. 
To leading order the critical moments of the family equal the critical
moments of the characteristic polynomials, up to a multiplicative
arithmetic constant.  For the family we average 
over the $L$-functions of conductor less
than~$X$, and then let~$X\to\infty$.  For the matrix groups
the averages are with respect to Haar measure. 
\end{enumerate}

The modeling described above will produce leading order asymptotics. 
To make more precise predictions requires other methods, which
are described in Section~\ref{sec:refined}.

\section{Elliptic curve families}\label{sec:EC}
For the remainder of the paper, $E$ is an elliptic curve over $\Q$,
with root number $w_E=\pm1$, and $L(s,E)$ is the $L$-function
associated to $E$ normalized so that $s=\frac12$ is the critical point.
If $F$ is any family of elliptic curves we write
$F=F^+\cup F^-$ where $F^+$ or $F^-$, respectively, are the
curves $E\in F$ with $w_E=+1$ or $w_E=-1$.
We write
$E_{a,b}$ for the curve $y^2=x^3+a x + b$.

Most of the information in this section can be found in recent papers
by Steven J.~Miller and Eduaro Due\~nez~\cite{M, MD1, MD2}, 
Nina Snaith~\cite{Sn1, Sn2} and Matthew Young~\cite{Y1, Y2, Y3}.  The author 
of this paper is just trying to convey the current understanding 
of the relationship between families of elliptic curve $L$-functions
and random matrix theory: he makes no claim
to any of the ideas presented here.

\subsection{Families with a given rank}

The following question is not well posed:

\emph{What is the correct random matrix model for the $L$-functions of a 
family of elliptic curves having a prescribed rank~$r$?}

The question is not well posed because there are (at least) 
two reasonable models, both of which seem to be appropriate for
certain families of elliptic curves.  We will first examine the
simplest case of rank~$r=1$.

Consider the following
families of rational elliptic curves:
\begin{equation}
\Fa(X)=\{E_{a,b}\ :\ |a|\le X^{\frac13},\ |b|\le X^{\frac12}\}
\end{equation}
and
\begin{equation}
\Fb(X)=\{E_{a,b^2}\ :\ |a|\le X^{\frac13},\ b^2\le X^{\frac12}\}.
\end{equation}
Note that the point $(0,b)$ on $E_{a,b^2}$ almost always has infinite
order, so almost all of the curves in $\Fb$ have rank
at least~1.

Let's consider $\Fa^-$ and $\Fb^-$.  
In both families we have $L(\frac12,E)=0$ because $w_E=-1$.
However, for $E\in \Fb^-$ we could have said
$L(\frac12,E)=0$ because $rank(E)\ge 1$.  The fact that the 
zero at $L(\frac12,E)$ for $E\in \Fb$ was \emph{constructed},
instead of just arising from parity considerations, has a
profound influence on the behavior of the $L$-function near
the critical point.

To understand the influence of the critical zero, we 
first consider the distribution of~$L'(\frac12,E)$,
which assuming standard conjectures is nonzero for 
almost all curves in our $w_E=-1$ families.
We have the following conjectures from~\cite{Y2}:
\begin{equation}\label{eqn:mattmom1}
\frac{1}{|\Fa^-(X)|} \sum_{E\in \Fa^-(X)} L'(\tfrac12,E)^k
\sim
c_1(k)(\log X)^{k(k+1)/2}
\end{equation}
while
\begin{equation}\label{eqn:mattmom2}
\frac{1}{|\Fb^-(X)|} \sum_{E\in \Fb^-(X)} L'(\tfrac12,E)^k
\sim
c_2(k)(\log X)^{k(k-1)/2}
\end{equation}

As the formulas show, the behavior at the critical point
is different for the two families, even though both families
could be described as ``a rank 1 family of elliptic curves.''
In particular, we see that the derivative $L'(\frac12,E)$
tends to be smaller for $E\in \Fb^-$.  
This can be explained by the tendency for the low-lying
zeros of $L(s,E)$ to be closer to the critical
point for~$E\in \Fb^-$.
That is, $\Fb^-$ should have more low-lying zeros, which
will cause the $L$-function to stay small near the critical point,
and so its derivative will also be small.
To make this idea precise we consider the 1-level density 
of the zeros.

Let $0<\gamma_{E,1}\le\gamma_{E,2}\le \gamma_{E,3}\le \cdots$ denote the 
imaginary parts of the zeros of $L(s,E)$ in the upper half of
the critical strip.  
Note that we have omitted the zero(s) at the critical point.
The \emph{one-level density} of the family $F(X)$ is 
defined to be the function $W_1$ which satisfies
\begin{equation}
\frac{1}{|F(X)|} \sum_{E\in F(X)} \sum_j \phi(\gamma_{E,j})
\sim
\int \phi(t)W_1(t) dt,
\end{equation}
as $X\to\infty$,
for nice functions~$\phi$.  That is, $W_1$ measures the density
of the zeros of the family.

\begin{figure}[ht]\label{fig:onelevel}
\begin{center}
\scalebox{0.7}[0.7]{\includegraphics{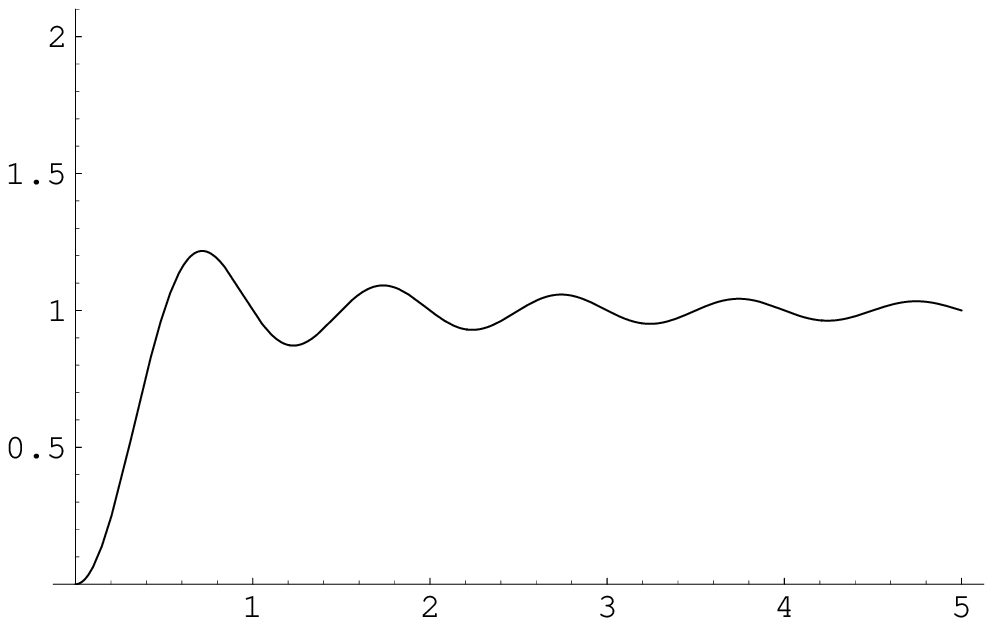}}
\hskip 0.4in
\scalebox{0.7}[0.7]{\includegraphics{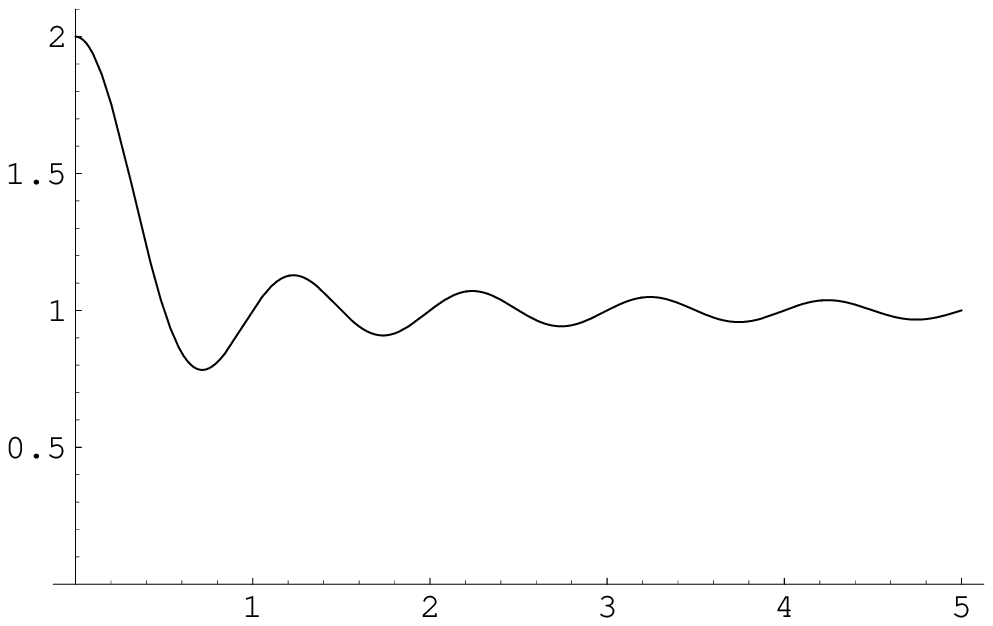}}
\\
\caption{\sf Conjectured one-level density of the 
noncritical $L$-function zeros of the family
$\Fa^-$ (left) and $\Fb^-$ (right).}
\end{center}
\end{figure}

The observation about the relative size of $L'(\frac12,E)$
can be restated as: the one-level density for the family
$\Fb^-$ should be more concentrated near 0 than the 
one-level density for the family~$\Fa^-$.  By using random matrix
theory and some other ideas we explain below, it is possible to 
produce a precise conjecture for the one-level densities of
these families.  These are given in Figure~\ref{fig:onelevel}.
The functions are rescaled so that the average spacing between
zeros is~$1$.  In the next section we explain where those 
conjectures came from.

\subsection{Two models for two kinds of families}\label{sec:twomodels}

The plots in Figure~\ref{fig:onelevel} are familiar.  The plot on
the left is the rescaled one-level density of the eigenvalues of matrices
from the group $SO(2N+1)$, in the limit as $N\to\infty$.
The plot on
the right is the rescaled one-level density of the eigenvalues of matrices
from the group $SO(2N)$, in the limit as $N\to\infty$.
Given those plot, what are our models for the two families?

Consider the following ways to make a polynomial $f(z)$ which has 
real coefficients, all its zeros on the unit circle, and 
(almost surely a simple) zero at $z=1$:
\begin{itemize}
\item{}The characteristic polynomial of a matrix in $SO(2N+1)$
\item{}$(z-1)$ times the characteristic polynomial of a matrix in $SO(2N)$
\end{itemize}

It should be clear that those two examples will have the one-level densities
pictured in Figure~\ref{fig:onelevel}.
These examples are the simplest cases of the two most commonly studied
higher rank families
of elliptic curve $L$-functions, which we now describe.

Suppose $E_T$ is a curve 
$y^2=x^3+a(T) x + b(T)$ of rank $r$ over $\Q(T)$.
Consider the following two families of rational elliptic curves:
\begin{equation}\label{eqn:allrankr}
\Fc(X)=\{E_{a,b} \ :\ |a|\le X^{\frac12},\ |b|\le X^{\frac13},\ \mathrm{rank}(E)\ge r\}
\end{equation}
and
\begin{equation}
\Fd(X)=\{E_{t}\ :\ |a(t)|\le X^{\frac12},\ |b(t)|\le X^{\frac13},\ t \in \N \}.
\end{equation}
As in our rank~1 example, we have $\Fd(X)\subset \Fc(X)$.  Again we consider
the subfamilies according to the sign of~$w_E$.
Let $F$ be either of the above rank~$r$ families.
If $r$ is odd then almost all the curves in~$F^-$ have rank $r$,
and almost all the curves in~$F^+$ have rank~$r+1$.  
If $r$ is even then the $F^+$ curves have rank~$r$ and the
$F^-$ curves have rank~$r+1$.
It is conjectured that $\Fd^+$ and $\Fd^-$ are approximately 
the same size provided $E$ 
has at least one place of multiplicative reduction. See~\cite{H}. 
If $r$ is even then it is possible that
$\Fc^+$ and $\Fc^-$ are approximately
the same size.
%And as in the rank~1 
%example it is believed that the two families have different behaviors 
%for their low-lying zeros.

We describe the two models which are believed to correspond to these families.
The names for these models was coined by Steven J.~Miller.

\subsection{Selecting to have zeros} \emph{The Interaction Model}

We are modeling a family that arises by restricting a much larger family
to a subfamily having at least~$r$ zeros at the critical point.
A matrix model for this family can be described as follows:
start with $SO(M)$ where $M=2N$ or $2N+1$ depending on whether we
are modeling $\Fc^+$ or $\Fc^-$, and
restrict to those matrices having $1$ as an eigenvalue of
multiplicity at least~$r$.
That is, you take matrices in $SO(M)$ and drag $r$ zeros to the 
critical point.

There are
some bad things about this model.  First, it is a set
of matrices, but it is not a group.  And while it is a perfectly
well-defined set, it is a measure zero subset of $SO(M)$, so
there is no canonical way to restrict Haar measure to it.

One solution, which has been analyzed by Snaith~\cite{Sn1} and
Due\~nez~\cite{D, MD1} is to first restrict to those matrices
which have $r$ eigenvalues in~$[-\varepsilon,\varepsilon]$,
and then let~$\varepsilon\to 0$.  The resulting measure is
the same as one obtains by taking Haar measure on $SO(M)$, 
formally substituting~$\theta=0$
for $r$ of the eigenvalues, and then omitting those terms which
vanish identically.  With the eigenvalues given by
$e^{i\theta_j}$, the induced measure on that set is
\begin{equation}
C(M, r)
\prod_{j=1}^M (1-\cos \theta_j)^r
\prod_{1\le j<k\le M} (\cos \theta_j-\cos \theta_k)^2
d\theta_1 \cdots d\theta_M
\end{equation}
where $C(M,r)$ is a normalization constant.

It is instructive to look at the one-level density for such matrices.
The one-level density is given by~\cite{Sn1, MD1}
\begin{equation}
\frac{\pi^2}{2}\theta \left(J_{r-\frac32}^2(\theta\pi)
+
J_{r-\frac12}^2(\theta\pi)
-
\frac{2r-1}{\theta\pi}J_{r-\frac12}(\theta\pi)J_{r-\frac32}(\theta\pi)\right) .
\end{equation}

\begin{figure}[ht]\label{fig:restrictonelevel}
\begin{center}
\scalebox{0.9}[0.9]{\includegraphics{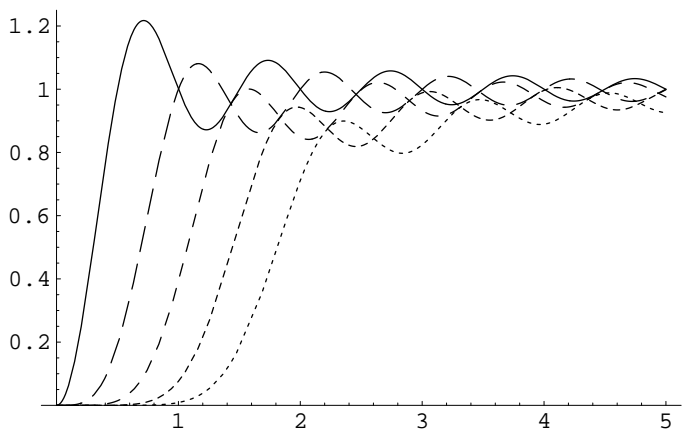}}
\\
\caption{\sf One-level density for $SO(M)$, restricted to have exactly
$r$ eigenvalues at~$\theta=0$,
for $r=1,2,3,4,5$.}
\end{center}
\end{figure}

There is some numerical evidence~\cite{MD1} that this model is accurate.
For this model it is possible to compute the critical moments of the
characteristic polynomials~\cite{Sn1}, but it does not seem that
all the ingredients are available
to use the recipe in~\cite{CFKRS} to conjecture the moments of the
family.

Note that when $r=0$ we recover~$SO(2N)$ and when $r=1$ we have
$SO(2N+1)$.  

\subsection{Imposing zeros} \emph{The Independent Model}

We are modeling a family that has an $r$th order zero at the
origin which arises from an explicit construction.  A model
for this situation can be found by assuming that the 
extra critical zeros are just inserted at the critical point, and all
the other zeros ignore them.
That is, start with a matrix in  $SO(M)$  where 
$M=2N-r$ or $2N+1-r$ depending on the parity of $r$ and the 
sign of~$w_E$.  Then the polynomial which models the $L$-function
is the characteristic polynomial of the matrix, multiplied
by~$(x-1)^r$.  The one-level density only depends on the
parity of $r$ and the sign of~$w_E$, and will be one of the 
functions shown in Figure~\ref{fig:onelevel}.

One can phrase the model strictly in terms of matrices by saying that
the model is given by the group
\begin{equation}
\begin{pmatrix}
I_{r\times r} & \\
 & SO(M) 
\end{pmatrix}
\end{equation}
where $I_{r\times r}$ is the $r\times r$ identity matrix and
$M=2N-r$ or $2N+1-r$ is chosen according to the parity of $r$ and
the sign of the functional equation.

There is numerical evidence~\cite{MD1}, level density calculations~\cite{M, MD1, Si, Y1}
and conjectures for moments~\cite{Y2} that this model gives accurate 
predictions for some specific families.

\subsection{Some issues}

It is worth repeating that the above models, even if they are correct,
are only intended to capture leading-term asymptotics.  Computer
experiments~\cite{MD1} find that the low lying zeros of the
family~$\Fd$ exhibit some anomalous behavior which presumably
will disappear when larger examples are computed.  It is 
entirely possible that family $\Fc$ is a union of families of
the form~$\Fd$, and this may contribute to a bias in those numerics.
In families of type~$\Fd$ the generators of the set of rational
points have very small height, and this may also introduce a bias. 
That is, the heights of the generators are on the order of the logarithm of
the conductor, while it is more typical to have the 
heights as large as a power of the conductor.
See Silverman~\cite{Si2}, Chapter 10, for a discussion of heights of generators.

It is possible that more accurate predictions of the one-level density
(using methods described in the next section) will show better agreement
with the data.  Those methods are also able to give extremely 
precise predictions for the moments of the $L$-functions~\cite{Y2},
and these give support to the models.

It is not clear that these families 
of elliptic curves give rise to families of $L$-functions 
as described in Section~\ref{sec:family}.
For the purpose of conjecturing the moments of the family
(which is why that definition of `family' was developed),
the key property is the orthogonality 
relation~\eqref{eqn:orthogonality}.
For specific families of type~$\Fd$ it should be possible
to evaluate such sums.  The result is likely to be quite
complicated, as in~\cite{Y2}.
A~subtle problem is that the parameter~$X$ in the elliptic
curve families is approximately the discriminant, not the
conductor.  By Szpiro's conjecture the
logarithm of the discriminant is within a factor of 6 of
the logarithm of the conductor, so it is possible that
ordering by discriminant is almost as good as ordering
by conductor.  It seems reasonable to model by setting~$N$, 
the size of the matrix,  equal to $\log X$, 
since what else would you choose?  If that choice is 
correct it suggests that $X$ is close to the
discriminant most of the time.  This has
been shown for some families~\cite{Y1}.

One can cook up an elliptic curve family which presumably is a hybrid of 
the models described here:  take a rank $r$ family of type $\Fd$ 
and restrict to those curves having rank at least $r+2$.   
If one makes the reasonable assumption that the ``extra'' zeros
created by this process do not interact with the original $r$ 
zeros imposed at the critical point, then one can use 
methods similar to~\cite{CKRS} to predict how many curves are
in the restricted family.  I am not advocating a reckless 
proliferation of elliptic curve models, but merely noting that
even if the two models described here are correct and can be 
refined to predict
lower order terms, they may not cover all families of interest. 

In Figure~\ref{fig:restrictonelevel} one can see that if $r$
is large then you are unlikely to find noncritical zeros close
to the critical point.  The name ``repulsion'' has
been given to this phenomenon.  The logic behind the name is that
the ``lowest zero'' is further from the critical point than
it would be if there were not a multiple critical zero.
Unfortunately, the ``lowest zero'' is not a well defined
object.  If you drag the lowest zero to the critical point then the
other zeros follow it toward the critical point, and at the moment
you increase the order of the critical zero there becomes a new
``lowest zero''.  If there is an $r$th order critical zero
and you count zeros correctly, then the ``lowest zero'' is actually
the $(r+1)$st zero, and it is likely to be \emph{closer} to the
critical point than a typical $(r+1)$st zero.  In other words,
the word ``attraction'' more accurately describes the situation!
There is no reason to change the current terminology, but keep in mind 
that in the model where one restricts to those matrices having
multiple eigenvalues at~$1$, the other eigenvalues have actually
moved \emph{towards} the critical point.

\section{Refined modeling}\label{sec:refined}

Random matrix theory is useful for making leading-order asymptotic 
predictions about
families of $L$-function.  
To understand the 
finer behavior of the family one must use heuristic techniques 
from number theory.

There are two main refinements to the leading term asymptotics.
First, $L$-function families are ordered by conductor, and 
we use the conductor to determine the appropriate
size matrices for our model.  There will be a discrepancy between
the limiting behavior for large matrices and the 
behavior for finite size matrices.
For some quantities, such as the nearest neighbor spacing of the 
zeros/eigenvalues, it is difficult to see the difference between the
asymptotics of the distribution  and the distribution for moderate
size matrices.  For other quantities, such as the value distribution,
for any computable range there is a notable difference between
the limiting quantity and the values that can be computed.  See
Keating and Snaith~\cite{KS1},  for a spectacular example
concerning the value distribution of $\Re\log\zeta(\frac12+it)$.
By a theorem of Selberg that is Gaussian 
in the limit $t\to \infty$, but for finite $t$ it
differs from a Gaussian in the same way as the characteristic 
polynomial of an appropriately sized random unitary matrix.

The second issue is the fact that there are 
lower order terms, and the lower order
terms for $L$-functions involve arithmetic factors, while the lower
order terms of random matrices do not.  Thus, the general shape of
expressions from random matrix theory can reveal what to expect
for $L$-functions, but the arithmetic ``correction terms'' must
be determined in some other manner. 
For zero spacings there are no arithmetic corrections in the leading
order terms.  For moments of $L$-functions the leading order
correction terms are fairly straightforward to determine~\cite{CG,CF,KS2}.
Just about everything else is quite subtle and one needs sophisticated
number-theoretic methods in order to make sensible conjectures.

For moments of $L$-functions
such conjectures are covered in detail in~\cite{CFKRS}.  
Matt Young~\cite{Y2} used these heuristics to compute the full main term for
various families of elliptic curve $L$-functions.  Those conjectures
give~\eqref{eqn:mattmom1} and \eqref{eqn:mattmom2} as special cases.  Thus, the
heuristics appear to correctly handle a variety of interesting
families.  (For the family $\Fc$ in \eqref{eqn:allrankr} 
our current understanding of the distribution of the coefficients~$a_p$
does not seem adequate to conjecture the moments of the family.
Even finding the leading order arithmetic factor seems difficult
in this case. However, random matrix calculations~\cite{Sn1} predict
the general shape of the moments.)

For quantities involving
zero statistics, moments are insufficient and one needs averages
of ratios of the $L$-functions.  This is addressed in~\cite{CFZ}.
For example, the expected value of the ratio
\begin{equation}\label{eqn:Lratio}
\frac{L(\frac12+\alpha,f)}{L(\frac12+\beta,f)}
\end{equation}
is sufficient to determine the one-level density of the family
$L(s,f)$, including the lower order correction terms due
to arithmetic effects.  See~\cite{CS} for many examples.
It should be possible to use these methods
to conjecture the ratio~\eqref{eqn:Lratio} for various higher rank
elliptic curve families, and thus give a precise conjecture
for the one-level density.  As of this writing this has not been done,
but probably it will have been done by the time this paper appears
in print.

In summary, by choosing the matrix group  and the size of the matrices 
appropriately, and with the appropriate arithmetic correction factor,
a family of $L$-functions can be modeled by the characteristic polynomials
of a collection of matrices.  The example of elliptic curve $L$-functions
where the elliptic curves are selected to have large rank shows that the
collection of matrices may not be a group.
In order to capture the lower order terms one must use heuristics
from number theory which do not explicitly involve random matrix theory.
Those heuristics also recover the leading order behavior which
previously required random matrix theory.

\end{document}